\documentclass[final]{siamltex}
\usepackage{amssymb,latexsym,amsmath}
\usepackage{mathrsfs}
\usepackage{epsfig}
\usepackage{caption}

\usepackage{dsfont}
\usepackage{color}
\usepackage{epsfig}
\usepackage{caption}

       \newtheorem{assumption}{\bf{Assumption}}[section]
\newtheorem{remark}{\bf{Remark}}[section]

\allowdisplaybreaks

\def\qed{\hfill $\diamond$}

\begin{document}


\title{On Stochastic Stability of a Class of non-Markovian Processes and Applications in Quantization}

\author{Serdar Y\"uksel
\thanks{Department of Mathematics and
    Statistics, Queen's University, Kingston, Ontario, Canada, K7L
    3N6.  Email: yuksel@mast.queensu.ca. This research was
    partially supported by the Natural Sciences and Engineering
    Research Council of Canada (NSERC).}
}

\maketitle

\begin{abstract}
In many applications, the common assumption that a driving noise process affecting a system is independent or Markovian may not be realistic, but the noise process may be assumed to be stationary. To study such problems, this paper investigates stochastic stability properties of a class of non-Markovian processes, where the existence of a stationary measure, asymptotic mean stationarity and ergodicity conditions are studied. Applications in feedback quantization and stochastic control are presented.
\end{abstract}

%



\section{Introduction}

Consider a stationary stochastic process $\{X_k, k \in \mathbb{Z}_+\}$ where each element $X_k$ takes values in some source space $\mathbb{X}$ (which we take to be $\mathbb{R}^n$ for some $n \in \mathbb{N}$ or some countable set) with process measure $\mu$, and a time-invariant update rule described by
\begin{eqnarray}\label{updateEq}
S_{k+1}=F(X_k, S_k)
\end{eqnarray}
where $S_k$ is an $\mathbb{S}$-valued state sequence (where we take $\mathbb{S}$ also to be $\mathbb{R}^n$ for some $n \in \mathbb{N}$ or some countable subset of $\mathbb{R}^n$), with $S_0=s$ or $S_0 \sim \kappa$ for some probability measure $\kappa$, independent of $X_k$. The question that we are interested in is whether for a given measurable and bounded $f$,
\begin{eqnarray}\label{convergence2}
\lim_{N \to \infty} {1 \over N} E[\sum_{k=0}^{N-1} f(X_{k},S_{k})]
\end{eqnarray}
or almost surely
\begin{eqnarray}\label{convergence3}
\lim_{N \to \infty} {1 \over N} \sum_{k=0}^{N-1} f(X_{k},S_{k})
\end{eqnarray}
exist and whether the limit is indifferent to the initial states/distributions. The function $f$ can be taken to be more general as follows:
\begin{eqnarray}\label{amstheorem}
\lim_{N \to \infty} {1 \over N} E[\sum_{k=0}^{N-1} f(X_{[k,\infty)},S_{[k,\infty)})]
\end{eqnarray}
or
\begin{eqnarray}\label{convergence4}
\lim_{N \to \infty} {1 \over N} E[\sum_{k=0}^{N-1} f(X_{(-\infty,k]},S_{k})]
\end{eqnarray}
Here, we use the notation that capital letters denote a random variable and small letters denote the realizations. We also have $y_{[m,n]}:=\{y_k, m\leq k \leq n\}$. One may also add another variable $U_k = g(S_k,X_k)$ where $U_k$ is an output of the system taking values in some set $\mathbb{U}$ and revise the formulation of the problem accordingly. We note that all of the random variables are defined on a common probability space $(\Omega, {\cal F}, P)$.

In (\ref{updateEq}) if $\{X_k\}$ were i.i.d, the process $\{S_k\}$ would be Markovian or if $\{X_k\}$ were Markovian, the joint process $\{(X_k,S_k)\}$ would be Markovian. For such Markov sources, there is an almost complete theory of the verification of stochastic stability through the analysis of finite-mean recurrence times to suitably defined sets (atoms, or artificial atoms constructed through {\it small} or {\it petite} sets and the splitting technique due to \cite{AthreyaNey} \cite{Nummelin}) as well as the regularity properties of the kernel (such as utilizing continuity of the transition kernel and majorization by a finite measure), see e.g. \cite{MeynBook} \cite{HernandezLermaLasserre}. For systems of the form (\ref{updateEq}) with only stationary $\{X_k\}$, however, there does not exist a complete theory even though the notion of {\it renovating events} \cite{Borovkov78} \cite{borovkov1992stochastically} that is related to the concepts of recurrence and coupling in Markov chains have been utilized in many applications especially in queuing theory. 

Such problems arise in many applications in feedback quantization and source coding, networks, and stochastic control. As an example, consider the following scheme which includes the $\Delta$-Modulation \cite{gersho1972stochastic} algorithm commonly used in source coding as a special case: Let $\{ X_k \}$ be stationary and ergodic, $Q: \mathbb{R} \to \mathbb{M} \subset \mathbb{R}, |\mathbb{M}| < \infty$ be a quantizer, and consider the following update:
\begin{eqnarray}
S_{k+1}=S_k + Q(X_k-S_k), \label{DeltaMod}
\end{eqnarray}
where $S_0=0$. Here, $S_k$ is the output of an adaptive encoder and $X_k$ is the source to be encoded. 

In addition to further adaptive coding schemes, applications include stabilization of controlled systems driven by noise processes with memory, design of networked control systems over channels with memory, as well as network and queuing systems. 

Such stability problems have been investigated for a number of setups; for an incomplete list see \cite{gersho1972stochastic,GoodmanGersho,Kieffer,Kieffer4,KiefferDunham,naraghi1990continuity,YukTAC2010,YukMeynTAC2010}. Notably, the contributions in Kieffer \cite{Kieffer} \cite{KiefferDunham} \cite{Kieffer4} are the most relevant ones to the discussion in this paper. These have studied problems motivated from applications in source coding and quantization as in (\ref{DeltaMod}). \cite{KiefferDunham} considered a non-Markovian setup where $\mathbb{S}$ is countable, \cite{Kieffer} considered a setup where $\mathbb{S}$ is not countable, but $f(x,\cdot)$ is continuous on $\mathbb{S}^{\mathbb{N}}$ for every $x$. Our approach and the proof technique is different than that considered in the literature; notably from that of Kieffer \cite{Kieffer}, and Kieffer and Dunham \cite{KiefferDunham} (as well as other contributions such as \cite{gersho1972stochastic} \cite{GoodmanGersho} \cite{YukTAC2010} and \cite{naraghi1990continuity} which can be approached by finite dimensional Markov chain formulations). 

Our approach builds on Markov process theoretic techniques where we model the stochastic process  $(X_{(-\infty,k]}, S_k)$ or $(X_{[k,\infty)}, S_{[k,\infty)})$ as an infinite dimensional Markov chain. The approach of viewing $(X_{(-\infty,k]},S_k)$ as a Markov chain, to our knowledge, first has been studied by Hairer \cite{hairer2009non}, where the focus of the author has been on the uniqueness of an invariant measure on the state process $S_k$, under the assumption that an invariant measure exists and further regularity assumptions. In this paper, we provide sufficient conditions for the existence of an invariant probability measure for the joint process while deriving our results. We also establish connections with asymptotic mean stationarity, in addition to the existence of an invariant measure, and ergodicity.

We will see that conditions of the form:
\begin{eqnarray}\label{tightS0}
\lim_{M \to \infty} \bigg(\limsup_{T \to \infty} {1 \over T} \sum_{k=0}^{T-1} P(|S_k| \geq M) \bigg)= 0,
\end{eqnarray}
play an important role for the stochastic stability results in this paper. Even though in the applications we consider we will explicitly study sufficient conditions for such a result, for a class of non-Markovian sources useful sufficient conditions (inspired from applications
in queuing and networks) are given in \cite{hajek1982hitting} and \cite{pemantle1999moment}. These
follow from Lyapunov-drift type conditions such as, with $S_n \geq 0$ for all $n$: $E[(S_{n+1} - S_n) 1_{S_n > L} | S_0, \cdots , S_n] < -A_1$, 
for some $A_1 >0, L < \infty$, and a uniform bound on jumps from below $L$ as $E[|S_{n+1}-S_n|^p | S_0,...S_n] < A_2$ 
for $p > 1$, $A_2 > 0$ leading to finite bounds on $\sup_n E[|S_n|^{p-1-\delta}]$ for arbitrarily small $\delta >0$, which through an application of Markov's inequality lead to (\ref{tightS0}). Thus, the findings in \cite{hajek1982hitting} and \cite{pemantle1999moment} together with the results in this paper can be used to obtain {\it Foster-Lyapunov} type drift criteria for various forms of stochastic stability.

A further related view to approach such problems is the traditional random dynamical systems view in which one studies the properties of the shifted sequences $(S_{[k,\infty)}, X_{[k,\infty)})$: Such a viewpoint leads to the interpretation that the entire uncertainty is realized in the initial state of the Markov chain, and the process evolves deterministically through a shift map. This approach has led to important contributions on ergodic theory and the introduction of useful notions such as asymptotic mean stationarity \cite{GrayProbabilit}. Connections between the two approaches and the implications on the convergences of (\ref{convergence2})-(\ref{convergence4}) will be made in the paper.

In Section \ref{sequenceMarkovI}, we discuss the conditions for the existence of an invariant probability measure. In Section \ref{ergodicityStat}, we discuss the conditions for asymptotic mean stationarity and ergodicity. This is followed by a study of applications in feedback quantization and networked stochastic control in Section \ref{Applications}. Section \ref{appendixA} in the Appendix contains a brief review of Markov chains and ergodic theorems for Markov chains.

\section{Stochastic stability of non-Markovian systems}\label{sequenceMarkovI}

Towards obtaining a method to study such systems, we will here view the process $(X_{(-\infty,k]}, S_k)$ as a $\mathbb{X}^{\mathbb{Z}_-} \times \mathbb{S}$-valued Markov process, similar to \cite{hairer2009non}. We recall that with $\mathbb{X}$ a complete, separable, metric (that is, a Polish) space, $\Sigma=\mathbb{X}^{\mathbb{Z}_-}$ is also a Polish space under the product topology. 

By a standard argument (e.g. {\sl Chapter 7} in \cite{durrett2010probability}), we can embed the one-sided stationary process $\{X_k, k \in \mathbb{Z}_+\}$ into a bilateral (double-sided) stationary process $\{ X_{k}, k \in \mathbb{Z} \}$. We first state the following.

\begin{lemma}
The sequence $(Z_k, S_k)$ with $Z_k = X_{(-\infty,k]}$ is a Markov process.
\end{lemma}
\textbf{Proof.} For any Borel $A \times B \in {\cal B}(\mathbb{X}^{\mathbb{Z}_-} \times \mathbb{S})$ and $k \geq 0$, the following holds almost surely:
\begin{eqnarray*}
&&P\bigg( (Z_{k+1},S_{k+1}) \in (A \times B) | Z_m, S_m, m \leq k \bigg) \\
&&= P\bigg( (Z_{k+1},F(X_k, S_k)) \in (A \times B) | Z_m, S_m, m \leq k \bigg) \\
&&= P\bigg( (Z_{k+1}\in A) \cap (F(X_k, S_k) \in B) | Z_m, S_m, m \leq k \bigg) \\
&&= P\bigg( X_{(-\infty,k+1]} \in A |Z_m, S_m, m \leq k\bigg) 1_{\{F(X_k, S_k) \in B \}} \\
&&= P\bigg( X_{(-\infty,k+1]} \in A |X_{(-\infty,m]}, S_m, m \leq k\bigg) 1_{\{F(X_k, S_k) \in B \}} \\
&&= P\bigg( X_{(-\infty,k+1]} \in A |X_{(-\infty,k]}\bigg) 1_{\{F(X_k, S_k) \in B \}} \\
&&= P\bigg( (X_{(-\infty,k+1]} \in A) \cap (F(X_k, S_k) \in B) | X_{(-\infty,k]}, S_k \bigg) \\
&&= P\bigg( (Z_{k+1},F(X_k, S_k)) \in (A \times B)  |Z_k,S_k \bigg)  \\
&&= P\bigg( (Z_{k+1},S_{k+1}) \in (A \times B) | Z_k, S_k \bigg)
\end{eqnarray*}
\qed

We let $\mathbb{P}$ denote the transition kernel for this process. Our inspiration for taking the approach below builds on the fact that, since $X_k$ is known to be stationary, if there were an invariant measure $v$ for this process, then this would decompose as
\[v(ds_0|x_{(-\infty,0]}) \pi(dx_{(-\infty,0]})\]
 with $\pi$ being the stationary measure for $X_k$. This follows by the invariance condition:
 \begin{eqnarray}
&&\int_{\mathbb{X}^{\mathbb{Z}_-} \times \mathbb{S}} \mathbb{P}(X_{(-\infty,k+1]},S_{k+1} \in B \times \mathbb{S} | x_{(-\infty,k]}, s_k) v(dx_{(-\infty,k]},ds_k) \nonumber \\
&&= \int_{\mathbb{X}^{\mathbb{Z}_-} \times \mathbb{S}} P(S_{k+1} \in \mathbb{S} | X_{(-\infty,k+1]} \in B, x_{(-\infty,k]}, s_k) \nonumber \\
&& \quad \quad \quad  \quad \quad \quad \times P(X_{(-\infty,k+1]} \in B | x_{(-\infty,k]}, s_k)  v(dx_{(-\infty,k]}, ds_k) \nonumber \\
&&=\int_{\mathbb{X}^{\mathbb{Z}_-} \times \mathbb{S}}  P(X_{(-\infty,k+1]} \in  B | x_{(-\infty,k]}, s_k)  v(dx_{(-\infty,k]},ds_k) \nonumber \\
&&=\int_{\mathbb{X}^{\mathbb{Z}_-} \times \mathbb{S}}  P(X_{(-\infty,k+1]} \in  B | x_{(-\infty,k]})  v(dx_{(-\infty,k]}, ds_k) \label{Invariance} \\
&&=\int_{\mathbb{X}^{\mathbb{Z}_-}}  P(X_{(-\infty,k+1]} \in  B | x_{(-\infty,k]})  v(dx_{(-\infty,k]}) \nonumber \\
&&= \pi(B) \label{invMarginal}
 \end{eqnarray}
 Here, in (\ref{Invariance}) we use the fact that given $x_{(-\infty,k]}$, to predict $X_{k+1}$, $S_k$ is non-informative. 

\subsection{Implications of the existence of an invariant probability measure}

If there is an invariant probability measure $\bar{P}$ for such a process we say that the process is {\it stochastically stable}. By the ergodic theorem (see Theorem \ref{convergenceT}), $\bar{P}$ almost surely
\begin{eqnarray}
\lim_{N \to \infty} {1 \over N} E_{x_{(-\infty,0]},s}[\sum_{k=0}^{N-1} f(X_{(-\infty,k]},S_{k})] = f^*(x_{(-\infty,0]},s) \label{asconv}
\end{eqnarray}
exists for all measurable and bounded $f$ and for corresponding functions $f^*$ (where the full set of convergence may depend on the function $f$).

The following assumption will be useful in establishing further stability results in Section \ref{ergodicityStat}. Recall that $S_0 \sim \kappa$ for some probability measure $\kappa$.
\begin{assumption}\label{absoluteCont}
The invariant measure $\bar{P}$ is such that $\pi \times \kappa \ll \bar{P}$. That is, $\bar{P}(A,B) =0$ implies that $\pi(A)\kappa(B) = 0$ for any Borel $A, B$. 
\end{assumption}

Under this assumption, we would have that the set of initial conditions which may not satisfy (\ref{asconv}) (this set has zero measure under $\bar{P}$) also has zero measure under the initial product probability measure $\pi \times \kappa$. Thus,
 \begin{eqnarray}
&& \int_{\mathbb{X}^{\mathbb{Z}_-}} \pi(dx) \int_{\mathbb{S}} \kappa(ds_0) E_{x_{(-\infty,0],s_0}} [{1 \over T} [\sum_{k=0}^{T-1} g(x_{(-\infty,k]},s_k)] \nonumber \\
&& \quad \quad \to \int_{\mathbb{X}^{\mathbb{Z}_-}} \int_{\mathbb{S}} \kappa(ds_0) \pi(dx) f^*(x_{(-\infty,0]},s_0)  \nonumber \\
 \end{eqnarray}

Furthermore, by Theorem \ref{convergenceT2}, sample paths also converge almost surely. Thus, convergence in the sense of (\ref{convergence2}), (\ref{convergence3}) and (\ref{convergence4}) will hold. We will later see that $(\ref{amstheorem})$ will also hold.

\subsection{Existence of an invariant probability measure with finite $\mathbb{S}$}

Our first result is for the setup with finite $\mathbb{S}$. For some related results and an alternative approach for the finite case, see \cite{Kieffer3}. 

\begin{theorem}\label{mainTheorem}
Consider the dynamical system given by (\ref{updateEq}). Suppose that $\mathbb{S}$ is finite. Then, the process is stochastically stable.
\end{theorem}

\textbf{Proof.}
Define for all $a \in \mathbb{S}$, the sequence of expected occupational measures
\begin{eqnarray}
&& v_t(dx_{(-\infty,k]} \times \{a\}) = E[{1 \over t} \sum_{k=0}^{t-1} 1_{\{X_{(-\infty,k]},S  \in dx_{(-\infty,k]} \times  \{a\} \}}] \nonumber \\
&& \quad \quad  = {1 \over t} \sum_{k=0}^{t-1} P(X_{(-\infty,k]}, S \in dx_{(-\infty,k]} \times \{a\}),\label{empiricmeasure}
\end{eqnarray}
where for every $k$, $X_{(-\infty,k]} \sim \pi$. Therefore, for any $t$, we can decompose $v_t(dx_{(-\infty,k]} \times \{a\} ) = \pi(dx_{(-\infty,k]}) v_t(a | x_{(-\infty,k]})$ since $\pi$ is a stationary measure by (\ref{invMarginal}). It follows then that \[v_t(dx_{(-\infty,k]} \times \{a\}) \leq \pi(dx_{(-\infty,k]})\] for every $a$ and $|\mathbb{S}| \pi(dx_{(-\infty,k]})$ is a majorizing {\it finite} measure for the sequence $v_t$. 

By \cite[Proposition 1.4.4]{HernandezLermaLasserre}, a sequence of probability measures which is uniformly countable additive is setwise sequentially precompact (see p. 6-8 in \cite{HernandezLermaLasserre}), a sufficient condition being that the sequence is majorized by a finite measure. Thus, $\{v_t\}$ has a converging subsequence $v_{t_k}$ so that for some probability measure $v$, $v_{t_k}(A) \to v(A)$ for all Borel $A$. Let $\mathbb{P}$ is the transition kernel for the Markov chain. Then through a Krylov-Bogoliubov-type argument, for every Borel $A$
\begin{eqnarray}
&& |v_{N}(A) - v_{N} \mathbb{P}(A)| \nonumber \\
&& = \bigg| {1 \over N} \bigg( (v_0(A)+\cdots+ v \mathbb{P}^{(N-1)}(A))  - (v_0 \mathbb{P}(A)+\cdots+v \mathbb{P}^{N}(A))  \bigg) \bigg| \nonumber\\
&& \leq {1 \over N} |v_0(A) - v_0 \mathbb{P}^{N}(A)|
 \to 0. \label{converEmp}
 \end{eqnarray}
Since $v_{t_k} \to v$, setwise, it follows that $v_{t_k}\mathbb{P}(B) \to v\mathbb{P}(B)$ also and hence $v(B)=v\mathbb{P}(B)$ and $v$ is stationary.  

As a result, the process has a converging subsequence setwise and the limit of this subsequence is invariant. Hence, there exists an invariant probability measure for the process. \qed

\subsection{Existence of an invariant probability measure with countable $\mathbb{S}$}

In this section, we assume that $\mathbb{S}$ is a countable set viewed as a subset of $\mathbb{R}$ whose elements are uniformly separated from each other; thus $\mathbb{S}$ is a uniformly discrete set in the sense that there exists $r > 0$ such that $|x-y| > r$ for all $x,y \in \mathbb{S}$ .

\begin{theorem}\label{mainTheoremII}
Consider the dynamical system given by (\ref{updateEq}). If (\ref{tightS0}) holds with the norm defined on $\mathbb{R}$, the process is stochastically stable.
\end{theorem}

\textbf{Proof.} 
(i) The sequence $v_t$ defined in (\ref{empiricmeasure}) is tight: As an individual probability measure, $\pi$ is tight. Since by (\ref{tightS0}), the sequence of marginals of $v_t^{\mathbb{S}}$ on $\mathbb{S}$ is also tight, it follows that the product measure is also tight: For every $\epsilon > 0$, there exists a compact set ${\cal L} \times {\cal M}$ in the product space so that 
\[ v_t\bigg(({\cal L} \times {\cal M})^C\bigg) \leq \pi({\cal L}^C) + v^{\mathbb{S}}_t({\cal M}^C) \leq \epsilon.\] This follows since $({\cal L} \times {\cal M})^C = ({\cal L}^C \times \mathbb{S}) \cup (\mathbb{X}^{\mathbb{Z}_-} \times {\cal M}^C)$, where for a set $A$, $A^C$ denotes its complement.

(ii) We show that the sequence $v_t$ is relatively compact under the $w$-$s$ topology \cite{balder2001, Schal}: Let $\mathbb{A}, \mathbb{B}$ be complete, separable, metric spaces. The $w$-$s$ topology on the set of probability measures ${\cal P}(\mathbb{A} \times \mathbb{B})$ is the coarsest topology under which $\int f(a,b) \nu(da,db): {\cal P}(\mathbb{A} \times \mathbb{B}) \to \mathbb{R}$ is continuous for every measurable and bounded $f$ which is continuous in $b \in \mathbb{B}$ for every $a \in \mathbb{A}$ (but unlike weak topology, $f$ does not need to be continuous in $a$).

Since the marginals on $x_{(-\infty,k]}$ is fixed, \cite[Theorem 3.10]{Schal} (see also \cite[Theorem 2.5]{balder2001}) establishes that the set of strategic measures is relatively compact under the $w$-$s$ topology when a tightness condition holds. By tightness from (i), let $v_{t_k}$ be a $w$-$s$ converging subsequence of $v_t$. Then, as in (\ref{converEmp}), for every Borel $A$
\begin{eqnarray}
&& |v_{N}(A) - v_{N} \mathbb{P}(A) | \nonumber \\
&& =  |{1 \over N} \bigg| (v_0(A)+\cdots+ v \mathbb{P}^{(N-1)}(A)) - (v_0 \mathbb{P}(A)+\cdots+v \mathbb{P}^{N}(A))  | \nonumber\\
&& \leq {1 \over N} (v_0(A) + v_0 \mathbb{P}^{N}(A)) \to 0, \label{con0}
 \end{eqnarray}

Now, $ v_{t_k}$ converges $w$-$s$ to $v$ for some $v$. In particular, for every measurable and bounded function $f$ (which is continuous in $s$ since $\mathbb{S}$ is countable), it holds that $\langle v_{t_k}, f \rangle \to \langle v, f \rangle$, where $\langle v_{t_k}, f \rangle := \int v_{t_k}(dx, s) f(x,s)$. We wish to show that $\langle v_{t_k} \mathbb{P},  f \rangle \to \langle v \mathbb{P} , f \rangle$, leading to the invariance of $v$ in view of (\ref{con0}). Now, let $f$ be measurable and bounded. We have that
\[ \int v_{t_k}(dx, s) f(x,s) \to \int v(dx, s) f(x,s)\]

Observe that the transitioned probability measure $v_{t_k}\mathbb{P}$ satisfies for every such $f$:
\begin{eqnarray}
 && \int v_{t_k}\mathbb{P}(dx, s) f(x,s) \nonumber \\
&& = \int \pi(dz) \bigg(\sum_{s'} v_{t_k}(s'|z) \int_x P(dx|z) f(x, F(z,s')) \bigg)
 \end{eqnarray}
 where $P(dx|z) = P(X_{(-\infty,k+1]} \in dx | X_{(-\infty,k]}=z)$. With $z=x_{(-\infty,k]}$, let \[\mathbb{P}f(z,s) = \int_x P(dx|z) f(x,F(x_k,s'))=:g(z,s).\]
We note here that with $z$ specified, $x_k$ is determined. The measurable function $g$ is continuous in $s$ for every $x$. This ensures that $\langle v_{t_k}, \mathbb{P} f \rangle = \langle v_{t_k} \mathbb{P},  f \rangle \to \langle v \mathbb{P} , f \rangle$ and for all bounded $f$ continuous in $s$: $\langle v, f \rangle = \langle v \mathbb{P} , f \rangle,$ and hence $v$ is invariant. \qed

\subsection{Existence of an invariant probability measure with $\mathbb{S}=\mathbb{R}^n$}
We have the following assumption.

\begin{assumption}\label{contin0'}
$F(x,s)$ is continuous in $s$ for every $x$. 
\end{assumption}

\begin{theorem}\label{mainTheoremIIb}
Consider the dynamical system given by (\ref{updateEq}). If (\ref{tightS0}) holds, under Assumption \ref{contin0'}, the process is stochastically stable.
\end{theorem}

\textbf{Proof.} 
(i) The sequence $v_t$ defined in (\ref{empiricmeasure}) is tight: As an individual probability measure, $\pi$ is tight. Since by (\ref{tightS0}), the sequence of marginals of $v_t^{\mathbb{S}}$ on $\mathbb{S}$ is also tight, it follows that the product measure is also tight: For every $\epsilon > 0$, there exists a compact set ${\cal L} \times {\cal M}$ in the product space so that 
\[ v_t\bigg(({\cal L} \times {\cal M})^C\bigg) \leq \pi({\cal L}^C) + v^{\mathbb{S}}_t({\cal M}^C) \leq \epsilon.\] This follows since $({\cal L} \times {\cal M})^C = ({\cal L}^C \times \mathbb{S}) \cup (\mathbb{X}^{\mathbb{Z}_-} \times {\cal M}^C)$, where for a set $A$, $A^C$ denotes its complement.

(ii) As before in Theorem \ref{mainTheoremII}, the sequence $v_t$ is relatively compact under the $w$-$s$ topology. 

Since the marginals on $x_{(-\infty,k]}$ is fixed, \cite[Theorem 3.10]{Schal} establishes that the set of strategic measures is relatively compact under the w-s topology under tightness. By tightness, let $v_{t_k}$ be a $w$-$s$ converging subsequence of $v_t$. Then, as in (\ref{converEmp}), for every Borel $A$
\begin{eqnarray}
&& |v_{N}(A) - v_{N} \mathbb{P}(A) | \nonumber \\
&& =  |{1 \over N} \bigg| (v_0(A)+\cdots+ v \mathbb{P}^{(N-1)}(A)) - (v_0 \mathbb{P}(A)+\cdots+v \mathbb{P}^{N}(A))  | \nonumber\\
&& \leq {1 \over N} (v_0(A) + v \mathbb{P}^{N}(A)) \to 0, \label{con0b}
 \end{eqnarray}

Now, $ v_{t_k}$ converges w-s to $v$ for some $v$. In particular, for every measurable and bounded function $f$ which is furthermore continuous in $s$ for every $x$, it holds that $\langle v_{t_k}, f \rangle \to \langle v, f \rangle$, where $\langle v_{t_k}, f \rangle := \int v_{t_k}(dx, ds) f(x,s)$. We wish to show that $\langle v_{t_k} \mathbb{P},  f \rangle \to \langle v \mathbb{P} , f \rangle$, leading to the invariance of $v$ in view of (\ref{con0b}). Now, let $f$ be measurable and bounded, but continuous in $s$. We have that
\[ \int v_{t_k}(dx, ds) f(x,s) \to \int v(dx, ds) f(x,s)\]

Observe that the transitioned probability measure $v_{t_k}\mathbb{P}$ satisfies for every measurable bounded $f$ continuous in $s$ for every $x$:
\begin{eqnarray}
 && \int v_{t_k}\mathbb{P}(dx, ds) f(x,s) \nonumber \\
&& = \int \pi(dz) \bigg(\int_{s'} v_{t_k}(ds'|z) \int_x P(dx|z) f(x, F(s',z)) \bigg)
 \end{eqnarray}
 where $P(dx|z) = P(X_{(-\infty,k+1]} \in dx | X_{(-\infty,k]}=z)$. With $z=x_{(-\infty,k]}$, let \[\mathbb{P}f(z,s) = \int_x P(dx|z) f(x,F(x_k,s'))=:g(z,s).\]
We note here that with $z$ specified, $x_k$ is determined. Since $F$ is continuous in $s'$ and $f$ is continuos is $s$, by the dominated convergence theorem $g$ is continuous in $s$ for every $z$. Thus, $\mathbb{P}$ preserves $w$-$s$ continuity and that $\langle v_{t_k}, \mathbb{P} f \rangle = \langle v_{t_k} \mathbb{P},  f \rangle \to \langle v \mathbb{P} , f \rangle$ and for all bounded $f$ continuous in $s$: $\langle v, f \rangle = \langle v \mathbb{P} , f \rangle,$ and hence $v$ is invariant. \qed

\subsection{Existence of an invariant probability measure under quasi-continuity conditions}

A large class of applications do not have the property that $\mathbb{S}$ is countable or that $F$ is continuous in $s$. To approach such problems in our framework, we impose the following quasi-Feller type condition which is natural for the applications we will consider.

\begin{assumption}\label{quasicontin}
$F(x,s)$ is continuous on $\mathbb{X} \times \mathbb{S} \setminus D$ where $D$ is a closed set with $P((X_{t+1},S_{t+1}) \in D | x_{(-\infty,t]}=x, s_t=s)= 0$ for all $x,s$. Furthermore, with $D_{\epsilon} = \{z: d(z,D) < \epsilon\}$ for $\epsilon > 0$ and $d$ the product metric on $\mathbb{X} \times \mathbb{S}$, for some $K < \infty$, we have that for all $x,s$ and $\epsilon > 0$
\[P\bigg( (X_{t+1},S_{t+1}) \in D_{\epsilon} | x_{(-\infty,t]}=x, s_t=s\bigg) \leq K \epsilon.\]
\end{assumption}
Note that we can write the above as
\[P\bigg( ((x_{(-\infty,t]},X_{t+1}),S_{t+1}) \in \{x_{(-\infty,t]}\} \times D_{\epsilon} | x_{(-\infty,t]}=x, s_t=s\bigg) \leq K \epsilon.\]
Furthermore, this is equivalent to the condition
\[P\bigg( ((x_{(-\infty,t]},X_{t+1}),S_{t+1}) \in \mathbb{X}^{\mathbb{Z}_-} \times D_{\epsilon} | x_{(-\infty,t]}=x, s_t=s\bigg) \leq K \epsilon.\]
We define $\mathbb{D}_{\epsilon} := \mathbb{X}^{\mathbb{Z}_-} \times D_{\epsilon}$. This is an open set in $\mathbb{X}^{\mathbb{Z}_-} \times \mathbb{S}$ and will be useful in the analysis to follow.

We remark that Assumption \ref{quasicontin} is related to what is referred to as the {\it quasi-Feller} condition introduced by Lasserre (see \cite[Section 7.3 ]{HernandezLermaLasserre}). Our definition here and the proof is different in part because we do not assume that the state space is locally compact. We have the following theorem.

\begin{assumption}\label{contin2'} 
\begin{itemize}
\item[(i)] If $\mathbb{X}$ is compact, $\int_{\mathbb{X}} P(X_{k+1} \in dx | X_{(-\infty,k]}=z) f(x)$ is continuous in $z$ for every continuous and bounded $f$ on $\mathbb{X}$.
\item[(ii)] If $\mathbb{X}$ is not compact, $\int_{\mathbb{X}} P(X_{k+1} \in dx | X_{(-\infty,k]}=z) f(x)$ is continuous in $z$ for every measurable and bounded $f$ on $\mathbb{X}$.
\end{itemize}
\end{assumption}

\begin{theorem}\label{mainTheorem3}
Suppose that Assumptions \ref{quasicontin} and \ref{contin2'} hold. If (\ref{tightS0}) holds, the system is stochastically stable.
\end{theorem}

\textbf{Proof.}
Assumption \ref{quasicontin} implies that every converging subsequence $v_{n_k}$ of
\[v_{n}(A \times B) = E_{x,s}[{1 \over n} \sum_{k=0}^{n-1} 1_{\{x_{(-\infty,k]},s_k \in (A \times B)\}}]\]
is such that for all $\epsilon > 0$
\[\limsup_{n_k \to \infty} v_{n_k}(\mathbb{D}_{\epsilon}) \leq K \epsilon.\]
Note that with $v=\lim_{n_k \to \infty} v_{n_k}$, it follows from the Portmanteau theorem (see e.g. \cite[Thm.11.1.1]{Dud02}) that
 \[v(\mathbb{D}_{\epsilon}) \leq K \epsilon.\]
Now, consider a weakly converging empirical occupation sequence $v_{t_k}$ and let this sequence have an accumulation point $v^*$. We will show that $v^*$ is invariant.

Observe that the transitioned probability measure $v_{t_k}\mathbb{P}$ satisfies for every continuous and bounded $f$:
\begin{eqnarray}
 && \int v_{t_k}\mathbb{P}(dx, ds) f(x,s) \nonumber \\
&& = \int \pi(dz) \bigg(\sum_{s'} v_{t_k}(ds'|z) \int_x P(dx|z) f(x, F(s',z)) \bigg)
 \end{eqnarray}
 where $P(dx|z) = P(X_{(-\infty,k+1]} \in dx | X_{(-\infty,k]}=z)$. With $z=x_{(-\infty,k]}$, let \[\mathbb{P}f(z,s) = \int_x P(dx|z) f(x,F(x_k,s'))=:g(z,s).\]
We note here that with $z$ specified, $x_k$ is determined. In the following, we argue that for continuous and bounded $f$, $g$ is continuous whenever $F$ is (thus outside $\mathbb{D}_{\epsilon}$).

If $\mathbb{X}$ is compact, by Tychonoff's Theorem $\mathbb{X}^{\mathbb{Z}_-}$ is locally compact. In this case, we will invoke \cite[Theorem 3.5]{serfozo1982convergence} for the following argument. If $z^n \to z$, \[f\bigg((z^n,x_{k+1}),F(x^n_k,s)\bigg) \to f\bigg((z,x_{k+1}),F(x_k,s)\bigg)\] for every $x_{k+1}$ and thus with $H_{s,n}(x_{k+1}):=f\bigg((z^n,x_{k+1}),F(z^n_k,s)\bigg)$, $H_{s}(x_{k+1}):=f\bigg((z,x_{k+1}),F(x_k,s)\bigg)$, it follows that $H_{s,n}(x^n_{k+1}) \to H_s(x_{k+1})$ as $z^n_{k+1} \to z_{k+1}$, thus we have {\it continuous convergence} as it is defined in \cite{serfozo1982convergence}. As a result, continuity of $g$ (outside $\mathbb{D}_{\epsilon}$) is established by a generalized dominated convergence theorem given in \cite[Theorem 3.5]{serfozo1982convergence} in view of weak continuity by Assumption \ref{contin2'}(i).  

If $\mathbb{X}$ is not compact, we invoke the generalized dominated convergence theorem of \cite[Theorem 2.4]{serfozo1982convergence}: Since $P$ is strongly continuous under Assumption \ref{contin2'}(ii), and $f$ is continuous and bounded, $g(z,s)=\int_x P(dx|z) f(x,F(x_k,s'))$ is continuous outside $\mathbb{D}_{\epsilon}$. 

Now, consider 
$\langle v_{t_k}, \mathbb{P}f \rangle = \langle v_{t_k}, g_f \rangle + \langle v_{t_k}, \mathbb{P}f - g_f \rangle$, where $g_f$ is a continuous function which is equal to $\mathbb{P}f$ outside an open neighborhood of $D$ and is continuous with $\|g_f\|_{\infty} = \|\mathbb{P}f\|_{\infty} \leq \|f\|_{\infty}$. The existence of such a function follows from the Tietze-Urysohn extension theorem \cite{Dud02}, where the closed set is given by $\mathbb{X}^{\mathbb{Z}_-} \times \mathbb{S} \setminus \mathbb{D}_{\epsilon}$. It then follows from Assumption \ref{quasicontin} that, for every $\epsilon > 0$ a corresponding $g_f$ can be found so that $\langle v_{t_k}, \mathbb{P}f - g_f \rangle \leq K \|f\|_{\infty} \epsilon$, and since $ \langle v_{t_k}, g_f \rangle \to  \langle v^*, g_f \rangle$, it follows that
\begin{eqnarray}
&& \limsup_{t_k \to \infty} | \langle v_{t_k}, \mathbb{P}f \rangle - \langle v^* , \mathbb{P}f \rangle | \nonumber \\
&& = \limsup_{t_k \to \infty} | \langle v_{t_k}, \mathbb{P}f - g_f \rangle - \langle v^* , \mathbb{P}f - g_f \rangle |  \nonumber \\
&& \leq \limsup_{t_k \to \infty} | \langle v_{t_k}, \mathbb{P}f - g_f \rangle| +  | \langle v^* , \mathbb{P}f - g_f \rangle | \nonumber \\
&& \leq 2 K' \epsilon
\end{eqnarray}
Here, $K' = 2 K \|f\|_{\infty}$ is fixed and $\epsilon$ may be made arbitrarily small. We conclude that $v^{*}$ is invariant. \qed

\begin{remark}
In his definition for quasi-Feller chains, Lasserre assumes the state space to be locally compact. A product space is locally compact if the individual coordinate spaces are compact by Tychonoff's Theorem and the state space in our formulation is not locally compact in general. However, we invoke tightness directly with no use of convergence properties of the set of functions which decay to zero as is done in \cite{HernandezLermaLasserre}. We also note that Gersho \cite{gersho1972stochastic} had obtained a similar result addressing points of discontinuity in the context of adaptive quantizer design.
\end{remark}

\section{Asymptotic mean stationarity and ergodicity}\label{ergodicityStat}

\subsection{Shifts and random dynamical systems view}
As an alternative approach, we may also view $X_{[k,\infty)},S_{[k,\infty)}$ as an infinite dimensional Markov chain. This viewpoint is more commonly adopted in the information theory literature (even though not explicitly stated as a Markov chain), as we discuss in the following. Such a view gives rise to important notions such as asymptotic mean stationarity. Note that such a viewpoint leads to the interpretation that the entire uncertainty is realized in the initial state of the Markov chain, and the process evolves deterministically. 

Let $\mathbb{X}$ be a complete, separable, metric space. Let ${\cal B}(\mathbb{X})$ denote the Borel sigma-field of subsets of $\mathbb{X}$, let $\Sigma=\mathbb{X}^{\mathbb{Z}_+}$ denote the sequence space of all one-sided (unilateral) infinite sequences drawn from $\mathbb{X}$. Thus, if $x \in \Sigma$ then $x = \{x_0,x_1, x_2, \dots\}$ with $x_i \in \mathbb{X}$. Let $X_n:\Sigma \to \mathbb{X}$ denote the coordinate function such that $X_n(x)=x_n$. Let $T$ denote the shift operation on $\Sigma$, that is $X_n(Tx)=x_{n+1}$. We also define the shift operation for a one-sided process defined on $\mathbb{Z}_-$ similarly: $X_{n-1}(Tx) = x_{n}, n \in \mathbb{Z}_-$.

With $\mathbb{X}$ a Polish space, $\Sigma=\mathbb{X}^{\mathbb{Z}_+}$ is also a Polish space under the product topology. Let ${\cal B}(\Sigma)$ denote the smallest $\sigma$-field containing all cylinder sets of the form $\{x: x_i \in B_i, m \leq i \leq n\}$ where $B_i \in {\cal B}(\mathbb{X})$, for all integers $m,n \geq 0$. Here, $\cap_{n \geq 0} T^{-n}{{\cal B}(\Sigma)}$ is the tail $\sigma-$field: $\cap_{n \geq 0} \sigma(X_n,X_{n+1},\cdots)$, since $T^{-n}(A) = \{x: T^nx \in A\}$.

Let $\mu$ be the measure on the process $\{X_0, X_1, \cdots\}$. This process is stationary and $\mu$ is said to be a stationary (or invariant) measure on $(\Sigma,{\cal B}(\Sigma))$ if $\mu(T^{-1}B)=\mu(B)$ for all $B \in {\cal B}(\Sigma)$. This random process is ergodic if $A=T^{-1}A$ implies that $\mu(A) \in \{0,1\}$.

\begin{definition} \cite{GrayKieffer} \label{AMSDefinition}
A process on a probability space $(\Omega, {\cal F}, {\bf P})$ with process measure $\mu$, is asymptotically mean stationary (AMS) if there exists a probability measure $\bar{P}$ such that
\begin{eqnarray}\label{amsdefn}
\lim_{N \to \infty} {1 \over N} \sum_{k=0}^{N-1} \mu(T^{-k} F) = \bar{P} (F),
\end{eqnarray}
for all events $F \in {\cal B}(\Sigma)$. Here $\bar{P}$ is called the stationary mean of $\mu$, and is a stationary measure.
\end{definition}

${\bar P}$ is stationary since, by definition $\bar{P}(F) = \bar{P}(T^{-1}F)$.

As elaborated on earlier, we may view $\{T^nX\}$ as a Markov chain (whose only uncertainty is hidden in the initial distribution) characterized by a transition function in the following, taking values in the Polish space $\mathbb{X}^{\mathbb{Z}_+}$. The kernel is given such that for every $n \in \mathbb{N}$
\[\int \mu(dx) P^n(x,B) = \mu(T^nX \in B) = \mu(X \in T^{-n}B).\]

We may define empirical occupation measures as follows:
\begin{eqnarray}\label{empiricalmeasure}
v_t(B) = {1 \over t} E_{\mu}[\sum_{k=0}^{t-1} 1_{\{T^kX \in B\}}] := {1 \over t} \sum_{k=0}^{t-1} \mu(X \in T^{-k}B)
\end{eqnarray}
In this case, the initial measure $\mu$ on the sequence space affects how convergence occurs. If the time averages converge {\it setwise} as in (\ref{amsdefn}) to some measure $\bar{P}$, the process is AMS.

It follows that if the system is AMS, (\ref{amstheorem}) holds since the set of simple functions is dense in the set of measurable and bounded functions under the supremum norm.

\subsection{Sufficient conditions for asymptotic mean stationarity}

It is an important question to ask when a process is AMS. 

\begin{theorem} \cite{rechard1956invariant} \label{Rechards} A process is AMS if and only if it is asymptotically dominated by a stationary process, that is there exists a stationary measure $\pi$ such that for Borel $B$ if $\pi(B)=0$ then $\lim_{n \to \infty} \mu(T^{-n}B)=0$.
\end{theorem}

Due to the Markov formulation, we can obtain the following direct condition to check whether the AMS property holds for systems of the form (\ref{updateEq}):

\begin{theorem}\label{amsimplied}
Let there exist a stationary measure $\bar{P}$ for the Markov chain $(X_{(-\infty,k]},S_k)$ for the system (\ref{updateEq}). Assumption \ref{absoluteCont} implies the AMS property for the process $(X_k, S_k)$. 
\end{theorem}

\textbf{Proof.} We will arrive at the conclusion using an ergodic theoretic result. From the last item of Theorem \ref{convergenceT3}, if $\pi \times \kappa \ll \bar{P}$, then the following uniform convergence holds:
\begin{eqnarray}
\lim_{T \to \infty} \sup_{f \in M(\mathbb{X}^{\mathbb{Z}_{-}} \times \mathbb{S}): \|f\|_{\infty} \leq 1} \bigg| E_{\pi \times \kappa} {1 \over T} [\sum_{k=0}^{T-1} f(X_{(-\infty,k]},S_k)]  -  \eta^*(f) \bigg|  = 0, \label{unifC}
\end{eqnarray}
for some invariant measure $\eta^*$ (not necessarily equal to $\bar{P}$), where $M(\mathbb{X}^{\mathbb{Z}_{-}} \times \mathbb{S})$ denotes the set of measurable and bounded functions on $\mathbb{X}^{\mathbb{Z}_{-}} \times \mathbb{S}$. 

We will now see that the above implies the AMS property. We will obtain the result for a general Markov process taking values in some Polish space $\mathbb{V}$ rather than the $\mathbb{X}^{\mathbb{Z}_-} \times \mathbb{S}$-valued process considered, for ease in presentation. 
Let the initial state be $v \in \mathbb{V}$ and the resulting measure on the state space $\mathbb{V}^{\mathbb{Z}_+}$ be $\mathbb{P}_{v}$. Let $\nu^*$ be an invariant probability measure for the Markov chain and $\mathbb{P}_{\nu^*}$ be the resulting stationary measure on the product space $\mathbb{V}^{\mathbb{Z}_+}$.  Consider a set $A \in {\cal B}(\mathbb{V}^{\mathbb{Z}_+})$ and let $A_1$ be a corresponding open finite dimensional set so that $x \in A$ is equivalent to $x_{[0,m]} \in A_1$ for some $m$. Then, 
\begin{eqnarray}
&& E_{\mathbb{P}_{v}} {1 \over T} [\sum_{k=0}^{T-1} 1_{\{ v_{k,\infty)} \in A\}}] = E_{\mathbb{P}_{v}} {1 \over T} [\sum_{k=0}^{T-1} 1_{\{ v_{k,k+m)} \in A_1\}}] \nonumber \\
&& = E_{\mathbb{P}_{v}} {1 \over T} \sum_{k=0}^{T-1} E[1_{\{ v_{k,k+m)} \in A_1\}} |v_{[0,k]}] = E_{\mathbb{P}_{v}} {1 \over T} [\sum_{k=0}^{T-1} E[1_{\{ v_{k,k+m)} \in A_1\}} |v_{k}]] \nonumber \\
&& = E_{\mathbb{P}_{v}} {1 \over T} [\sum_{k=0}^{T-1} g(v_k)] = E_{v_0=v} {1 \over T} [\sum_{k=0}^{T-1} g(v_k)] \nonumber \\
&& \to \int \eta^*(dv)g^*(v) = \int \eta^*(dv) g(v) \nonumber \\
&& = E_{\mathbb{P}_{\eta^*}} {1 \over T} [\sum_{k=0}^{T-1} 1_{\{ v_{k,\infty)} \in A\}}] = \mathbb{P}_{\eta^*}(A),
 \end{eqnarray}
 for some invariant measure $\eta^*$ and measurable $g^*$. Here, $g(v_k)=E[1_{\{ v_{k,k+m)} \in A_1\}} |v_{k}]$. The first equality above follows from the fact that $A$ is a finite-dimensional cylinder set and $A_1$ is the corresponding finite dimensional set, the second equality from the iterated expectations, the third from the fact that $v_k$ is Markov. Since, $g(v)$ is measurable and bounded, (\ref{unifC}) leads to the desired result. 
 
 \qed

 Naturally, if $\bar{P}$ is the unique invariant measure, $\eta^* = \bar{P}$ in (\ref{unifC}). Next, we study the uniqueness problem.

\subsection{Ergodicity}

Ergodicity is a desirable stability property, since it allows for the sample path averages to converge to the same limit in the ergodic theorem regardless of the initial distribution, leading to crucial consequences in information theoretic and control theoretic applications. For a Markov chain, the uniqueness of an invariant probability measure implies ergodicity (see e.g. \cite[Chp.~2]{HernandezLermaLasserre}). With the random dynamical systems view, sufficient conditions such as mixing can be utilized \cite{GrayProbabilit}, however these may be restrictive. 

Consider an $\mathbb{X}$-valued Markov chain with transition kernel $P$, where $\mathbb{X}$ is a complete, separable and metric space. 
\begin{definition}
A Markov chain is $\mu$-irreducible, if for any set $B \in {\cal B}(\mathbb{X})$ such that $\mu(B)>0$, and $\forall x \in \mathbb{X}$, there exists some integer $n>0$, possibly depending on $B$ and $x$, such
that $P^n(x,B) > 0$, where $P^n(x,B)$ is the transition probability in $n$ stages from $x$ to $B$.
\end{definition}

A maximal irreducibility measure $\psi$ is an irreducibility measure such that for all other irreducibility measures $\phi$, we have $\psi(B)=0 \Rightarrow \phi(B)=0$ for any $B \in \mathcal{B(X)}$ (that is, all other irreducibility measures are absolutely continuous with respect to $\psi$). Whenever a chain is said to be irreducible, irreducibility with respect to a maximal irreducibility measure is implied. A maximal irreducibility measure $\psi$ exists for a $\mu$-irreducible Markov chain, see \cite[Propostion 4.2.4]{MeynBook}. The following is a well-known result.

\begin{theorem}
Let $\{X_t\}$ be a $\psi$-irreducible Markov chain which admits an invariant probability measure. The invariant measure is unique.
\end{theorem}

\textbf{Proof.} 
Let there be two invariant probability measures $\mu_1$ and $\mu_2$. Then, there exists two {\it mutually singular} invariant probability measures $\nu_1$ and $\nu_2$, that is $\nu_1(B_1)=1$  and $\nu_2(B_2)=1$, $B_1 \cap B_2 = \emptyset$ and that $P^n(x,B_1^C)=0$ for all $x \in B_1$ and $n \in \mathbb{Z}_+$ and likewise $P^n(z,B_1^C)=0$ for all $z \in B_1$ and $n \in \mathbb{Z}_+$ (see e.g. \cite[Lemma 2.2.3]{HernandezLermaLasserre}). This then implies that the irreducibility measure has zero support on $B_1^C$ and zero support on $B_2^C$ and thus on $\mathbb{X}$, leading to a contradiction. 
\qed

A complementary condition for ergodicity is the following.

\begin{definition}
For a Markov chain with transition kernel $P$, a point $x$ is accessible if for every $y$ and every open neighbourhood $O$ of $x$, there exists $k > 0$ such that $P^k(y,O) > 0$.
\end{definition}

One can show that if a point is accessible, it belongs to the (topological) support of every invariant measure (see, e.g., Lemma 2.2 in \cite{Hairer}). Recall that the support (or spectrum) of a probability measure is defined to be the set of all points $x$ for which every open neighbourhood of $x$ has positive measure.

We recall that a Markov chain $V_t$ is said to have the strong Feller property if $E[f(V_{t+1})|V_t=v]$ is continuous in $v$ for every measurable and bounded $f$. 

\begin{theorem}\label{uniquenessInv} \cite{Hairer} \cite{da1996ergodicity}
If a Markov chain over a Polish space has the strong Feller property, and if there exists an accessible point, then the chain can have at most one invariant probability measure.
\end{theorem}

However, a Markov chain defined as $(X_{(-\infty,k]},S_k)$ cannot be strongly Feller due to the memory in the source: Take $f(x)=1_{\{x_{-1} \in A\}}$ (where $x=x_{(-\infty,0]}$) for some closed set $A$, then $E[f(X_{(-\infty,k+1]})|x_{(-\infty,k]}] = 1_{\{x_k \in A\}}$ is not continuous. 

Nonetheless, we can have the following slight generalization.
\begin{theorem}\label{strongC}
Suppose that $E[f(X_{k+1},S_{k+1})|X_{(-\infty,k]} = x_{(-\infty,k]},S_k=s_k]$, for measurable and bounded $f: \mathbb{X} \times \mathbb{S} \to \mathbb{R}$, is continuous in $(x_{(-\infty,k]}, s_k)$. Suppose further that there exists an accessible point for the Markov chain $\{(X_{(-\infty,k]},S_k)\}$. The chain can have at most one invariant probability measure. 
\end{theorem}

\textbf{Proof.} As in \cite{Hairer}, suppose there exist two different invariant probability measures $\pi_1, \pi_2$ both of which must include $(x,s)$ in their topological supports. Then, there exist disjoint sets $U$ and $V$ and probability measures $\tilde{\pi}_1$ and $\tilde{\pi}_2$ so that $\tilde{\pi}_1(U)=1$ and $\tilde{\pi}_2(V)=1$ (see e.g. \cite[Lemma 2.2.3]{HernandezLermaLasserre}). Now, there cannot exist a Borel $A \subset \mathbb{X} \times \mathbb{S}$ such that under one measure it puts $P((x_{(-\infty,t]},X_{t+1},S_{t+1}) \in (\{x_{(-\infty,t]}\} \times A)|x_{(-\infty,t]},s)  = P((X_{t+1},S_{t+1}) \in A) |x_{(-\infty,t]},s)=1$ $\tilde{\pi}_1$ a.s. and $P((X_{t+1},S_{t+1}) \in A)|x_{(-\infty,t]},s) = 0$ $\tilde{\pi}_2$ a.s since the function $E[1_{\{(X_{t+1},S_{t+1}) \in A\}} |x,s]$ is continuous in $x,s$ and if a continuous function is a constant $\tilde{\pi}_i$ almost everywhere, then it should be a constant in the topological support of the probability measure. The conditions of the theorem imply that $\mathbb{S}$ is a uniformly separated set. By an iterated analysis it follows that for every finite dimensional cylinder set on $\{X_k, S_k\}$, the supports of the measures induced under $\tilde{\pi}_1$ and $\tilde{\pi}_2$ on these finite dimensions must be consistent; and by stationarity of $\{X_t\}$, equal. This implies that the measures $\tilde{\pi}_1$ and $\tilde{\pi}_2$ must be equal. 
\qed

For applications such as $\Delta$-Modulation, however, we will see that the continuity assumption in Theorem \ref{strongC} fails to hold. To be able to apply the result for such setups, we have the following relaxation utilizing Assumption \ref{quasicontin}.

\begin{theorem}\label{strongC2}
Suppose that for measurable and bounded $f: \mathbb{X} \times \mathbb{S} \to \mathbb{R}$, $E[f(X_{k+1},S_{k+1})|X_{(-\infty,k]} = x_{(-\infty,k]},S_k=s_k]$ is continuous  in $(x_{(-\infty,k]},s_k)$ for all $(x,s) \in (\mathbb{X}^{\mathbb{Z}_-} \times \mathbb{S}) \setminus D$ for some closed set $D$ which satisfies the conditions in Assumption \ref{quasicontin}. Suppose further that there exists an accessible point $(x,s) \notin D$ for the Markov chain $\{(X_{(-\infty,k]},S_k)\}$. The chain can have at most one invariant probability measure. 
\end{theorem}

\textbf{Proof.} The proof follows from that of Theorem \ref{strongC}, despite the presence of the discontinuity set $D$. \qed 

\section{Applications}\label{Applications}

In this section, we consider applications in feedback quantization and networked control. 

\subsection{Adaptive Quantization}

Adaptive quantization for stationary sources has been studied in particular in \cite{Kieffer}, \cite{KiefferDunham} and \cite{gersho1972stochastic}. This paper generalizes the results of \cite{gersho1972stochastic} which investigated $\Delta$-Modulation only for finite order Markov sources. We believe that the approach in this paper is more accessible than the arguments in \cite{Kieffer} and \cite{KiefferDunham} in part because it allows for, through a unified approach, a Markov chain theoretic approach and also leads to an ergodicity analysis in addition to asymptotic mean stationary.

\subsubsection{$\Delta$-Modulation} \hfill
\begin{theorem}\label{DeltaMAMS}
Let $X_k$ be stationary and ergodic $\mathbb{R}$-valued process stationary process measure $\pi$, $Q: \mathbb{R} \to \{-m, m\}$, with the following update:
\[S_{k+1}=S_k + Q(X_k-S_k),\]
where $S_0=0$ and $Q(Z)= m 1_{\{Z \geq 0\}} -m 1_{\{Z < 0\}}$. Suppose further that $E[Q(X_0 - m)] < 0$ and $E[Q(X_0+m)] > 0$ (equivalently $P(X_0 \geq m) < 1/2, P(X_0 \leq -m) < 1/2$). Then, the system is stochastically stable in the sense that there exists an invariant probability measure. Furthermore, if for every $m,k$, and non-empty open $A_k$, $\pi(X_{[m,k]} \in \prod_{t=m}^{k} A_k) > 0$, the system is AMS. If in addition $E[g(X_1)| x_{(-\infty,0]}]$ is continuous in $x_{(-\infty,0]}$ for measurable and bounded $g$, $(X_k,S_k)$ is ergodic.
\end{theorem}

Note that here $\mathbb{S} = \{km, k \in \mathbb{Z}\}$ is a countable set.

An example where the {\it Lebesgue-irreducibility} type condition ($\pi(X_{[m,k]} \in \prod_{t=m}^{k} A_k) > 0$) holds is
\[X_{t+1}=\sum_{i=0}^{\infty} \alpha_i W_{t-i},\]
with $\sum_t |\alpha^2_t| < \infty$ and $W_t$ is a sequence of i.i.d. Gaussian random variables. 

An example where the continuity condition for ergodicity holds is the following auto-regressive representation
\[X_{t+1}=\sum_{i=0}^{N-1} \alpha_i X_{t-i} + W_t,\]
with the roots of $1 - \sum_{i=1}^N \alpha_{i-1} z^{-i}$ strictly inside the unit circle and $W_t$ a sequence of i.i.d. Gaussian random variables. This follows since 
\[E[g(X_1)| x_{(-\infty,0]}] = \int g(z) \eta(z- \sum_{i=0}^{N-1} \alpha_i x_{-i}) dz,\]
with $\eta$ denoting the Gaussian density and by an application of the dominated convergence theorem, this expression is continuous in $x_{(-\infty,0]}$.

\textbf{Proof.}
Observe first that 
\begin{eqnarray}\label{tightS01}
&& \lim_{M \to \infty} \bigg(\limsup_{T \to \infty} {1 \over T} \sum_{k=0}^{T-1} P(|S_k| \geq M) \bigg) \nonumber \\
&& = \lim_{M \to \infty} \bigg(\limsup_{T \to \infty} {1 \over T} \sum_{k=0}^{T-1}  (1 - P(|S_k| < M)) \bigg) \nonumber \\
&& = \lim_{M \to \infty} \bigg(1 - \bigg(\liminf_{T \to \infty} {1 \over T} \sum_{k=0}^{T-1}  P(|S_k| < M) \bigg) \bigg) \nonumber \\
&& = 1 - \lim_{M \to \infty} \bigg(\liminf_{T \to \infty} {1 \over T} \sum_{k=0}^{T-1}  P(|S_k| < M) \bigg) \nonumber 
\end{eqnarray}
and thus (\ref{tightS0}) can be equivalently written as
\begin{eqnarray}\label{tightS02}
\lim_{M \to \infty} \bigg(\liminf_{T \to \infty} {1 \over T} \sum_{k=0}^{T-1}  P(|S_k| < M) \bigg) = 1.
\end{eqnarray}

Using the fact that by the ergodicity of the source the following hold almost surely:
\begin{eqnarray}
&&\lim_{n \to \infty} \bigg( m | \{- n \leq -i : |X_i| \leq m\}|  - m  | \{- n \leq -i : |X_i| > m\}| \bigg) =  \infty, \nonumber \\
&& \lim_{n \to \infty} \bigg( -m | \{- n \leq -i : |X_i| \leq m\}|  + m  | \{- n \leq -i : |X_i| > m\}| \bigg) =  -\infty, \nonumber 
\end{eqnarray}
Kieffer and Dunham's \cite[Theorem 2]{KiefferDunham} shows that the coding scheme satisfies the condition \cite[Eqn. (2.2)]{KiefferDunham}, which in turn implies (\ref{tightS02}). This follows since, with $K$ a finite set, the condition $S_i \in K$  for some $i \in \{n, \cdots, n+N\}$ \cite[Eqn. (2.2)]{KiefferDunham} implies that $|S_i| \leq K_1 + Nm \leq K_2 N$ for constants $K_1, K_2$ since $|S_i - S_j| \leq |i-j| m$. By Theorem \ref{mainTheoremII}, the system is stochastically stable.

{\it Asymptotic mean stationarity:} For the AMS property, we show that Assumption \ref{absoluteCont} holds: Let $X_{[-m, 0]} \in B$ and $S_0=0$ have a zero measure under $\bar{P}$. Then, $\pi(B) = 0$. To show this, consider the contrapositive: If $\pi(B) > 0$, by the condition that all finite-dimensional cylinder sets consisting of non-empty open sets have positive measure conditioned on any past event, it follows that for some $S_0=s^*$ with positive measure under $\bar{P}$, there exists a positive probability event $X_{[0,m]} \in B$ so that $S_m=0$. With, 
\begin{eqnarray}
&& \bar{P}(X_{[-m,0]} \in B, S_0=0) = \bar{P}(X_{[0,m]} \in B, S_m=0)  \nonumber \\
&& \quad \quad \quad \quad \geq \int_{z} \bar{P}(dz,s^*) \mathbb{P}(x_{[0,m]} \in B, S_m=0| z,s^*) > 0,
\end{eqnarray}
it follows that the absolute continuity condition holds, and by Theorem \ref{amsimplied}, the AMS property.

{\it Ergodicity}: We can establish the uniqueness of an invariant probability measure through either irreducibility properties or the following argument. Consider the point $p_0=\{m/2\}^{\mathbb{Z}_-} \times \{0\}$. We argue that this point is accessible. Recall that an open set in a product topology is a Borel set in the product space consisting of finitely many open sets with the rest being $\mathbb{X}$ itself or arbitrary union of such sets. Now, consider any $x_{(-\infty,0]},s$. From this point, we will show that for every open neighborhood $U$ of $p_0$, there exists some $k > 0$ so that $P(X_{(-\infty,k]}  \in U|x_{(-\infty,0]}, s) > 0$. For $x \in U$ for such $U$, $x_{(-\infty,0]} \in \prod_{l=-\infty}^0 A_l$ for finitely many non-empty open sets which are not equal to $\mathbb{X}$ and the rest being $\mathbb{X}$ (see e.g. \cite[Theorem 2.4.4]{Dud02}). Let $-l$ be the largest index for which $A_{-l} \neq \mathbb{X}$. Hence, it is evident that $x_{(-\infty,l]}$ can take values in this open set for a given $x_{(-\infty,0]}$. We also need to ensure that $S_l$ hits zero. To allow for this to happen, we further shift the process to the left: for any sufficiently small $\epsilon >0$, identify a sequence of events from $r$ to $r+l$ so that $S_{r+l}=0$ for some $S_r=s$ when $|X_k - m/2| \leq \epsilon$ in this time interval. As a result, at time $r+l$ the state process hits $0$ and the process $X_{(-\infty,r+l]}$ hits the open set with positive probability with $S_{r+l}=0$. Finally,  continuity holds due to the continuity of the noise process: The sets of points where continuity fails, $D=\{x: x = km, k \in \mathbb{Z}\}$, is a closed set satisfying the conditions in Theorem \ref{strongC2}, and $p_0$ is outside this set. By Theorem \ref{strongC2}, the process is ergodic. \qed

\subsubsection{Adaptive Quantization of Goodman and Gersho}

Consider the following update equations \cite{GoodmanGersho}:
\begin{eqnarray}
&& V_t = \Delta_t Q_1(X_t/\Delta_t) \nonumber \\
&& \Delta_{t+1}=\Delta_t Q_2({|X_t| \over \Delta_t}), \quad \Delta_0 = b \label{updateGG}
\end{eqnarray}
Here, $\Delta_t$ is the bin size of the uniform quantizer with a finite range and $|Q_1(\mathbb{R})| < \infty, |Q_2(\mathbb{R}_+)|< \infty$. $V_t$ is the output which is to {\it track} the source process $X_t$. Suppose further that $Q_2$ is non-decreasing.

\begin{theorem}
Let $X_t$ be a stationary and ergodic (non-deterministic) Gaussian sequence, $\zeta=\lim_{x \to \infty} Q_2(x) > 1$, $Q_2(0) = \lim_{x \downarrow 0} Q_2(x)< 1$ and $\log_2(Q_2(\cdot)) \in \mathbb{Q}$. Then, the system is stochastically stable. If in addition, with $\{\alpha_1, \alpha_2, \cdots, \alpha_L\}$ a set of pairwise relatively prime integers and $\log_2(Q_2(\cdot)) \in \{\alpha_k m\}$ for some $m \in \mathbb{Q}$, the process is AMS, and furthermore, ergodic. 
\end{theorem}

\textbf{Proof.} Consider 
\[\log_2(\Delta_{t+1}) = \log_2(\Delta_t) + \log_2(Q_2({|X_t| \over \Delta_t}))\]
$\log_2(\Delta_t)  - \log_2(\Delta_0) \in \mathbb{Q}$ for all $t$. Let $S_t = \log_2(\Delta_t)$. This sequence takes values in a countable set and satisfies
\[ S_{t+n} - S_{t} = \sum_{k=t}^{t + n-1} \log_2(Q_2({|X_k| \over \Delta_k})).\]

As in the proof of Theorem \ref{DeltaMAMS}, \cite[Theorem 4]{KiefferDunham} shows that the coding scheme satisfies the condition \cite[Eqn. (2.2)]{KiefferDunham}, which implies (\ref{tightS02}). By Theorem \ref{mainTheoremII}, the system is stochastically stable.



{\it The AMS property}: Since $\{\alpha_k\}$ is a set of numbers that are relatively prime $\mathbb{S}$ consists of all integer multiples of $m$ shifted by the initial value $\log_2(b)$. This follows from the property of relatively prime numbers due to B\'ezout's lemma; see \cite[Lemma 7.6.2]{YukselBasarBook}. The argument for the AMS property then follows as before through the absolute continuity condition: Any invariant measure is such that $\bar{P}(\cdot,s) \ll \bar{P}(\cdot,s')$ for all admissible $s, s'$ and by Theorem \ref{amsimplied}, the result follows.

{\it Ergodicity}: In this case, the point $(\{0\}^{\mathbb{Z}_-},\log_2(b))$ is accessible by the same arguments adopted in the proof of Theorem \ref{DeltaMAMS} and the steps leading to the AMS property above. By Theorem \ref{strongC2}, the process is ergodic. \qed

\subsection{Stochastic networked control}

We consider a stabilization problem in stochastic networked control where a linear system is controlled over a communication channel. We will study the approach in \cite{YukTAC2010}, \cite{YukMeynTAC2010} (see \cite{YukselBasarBook} for a detailed discussion). Consider the following control system, with $U_t$ a control variable,
\begin{eqnarray}\label{scalarSystem}
X_{t+1}=aX_t + bU_t + W_t.
\end{eqnarray}
where $|a| \geq 1$, $W_t$ is i.i.d, admitting a probability measure $v$ which admits a density, positive everywhere and bounded. Furthermore, $E[|W_t|^{2+\zeta}] < \infty$ for some $\zeta > 0$.

In the application considered, a controller has access to quantized information from the state process. The quantization is described as follows. An adaptive quantizer has the following form with $Q_K^{\Delta}$  being a uniform quantizer with $K+1$ bins and bin-size $\Delta$, $Q^{\Delta}_K: \mathbb{R} \to \mathbb{R}$ satisfies the following for $k=1,2\dots,K$:
\begin{eqnarray}
Q_K^{\Delta}(x) = \begin{cases}   (k - {1 \over 2} (K + 1) ) \Delta, \quad \quad  \mbox{if} \quad x \in [ (k-1-{1 \over 2} K   ) \Delta , (k-{1 \over 2} K  ) \Delta)  \\
 {1 \over 2} (K - 1) \Delta,  \quad \quad \mbox{if} \ \ x = {1 \over 2} K \Delta \\
  0 ,  \quad \quad
  	\mbox{if} \ \ x \not\in [- {1 \over 2} K  \Delta, {1 \over 2} K  \Delta]
\end{cases} \nonumber
\end{eqnarray}
%
%
With $K = \lceil |a|+\epsilon \rceil$, $R = \log_2(K + 1)$, let $R'=\log_2(K)$. We will consider the following coding and quantization update policy. For $t \geq 0$ and with $\Delta_0 > L$ for some $L \in \mathbb{R}_+$, and $\hat{x}_0 \in \mathbb{R}$, consider:
\begin{eqnarray}\label{QuantizerUpdate2CHP7}
&&U_t = - {a \over b} \hat{X}_{t}, \quad \hat{X}_t = Q_K^{\Delta_t}(X_t), \nonumber \\
&& \quad \Delta_{t+1} = \Delta_t \bar{Q}(|{ X_{t} \over \Delta_t 2^{R'-1} }|, \Delta_t) \nonumber
\end{eqnarray}
Suppose that with $\delta, \epsilon, \alpha > 0$ with $\alpha < 1$ and $L > 0$ the following hold
\begin{eqnarray}\label{QuantizerNoiselessUpdateEquation}
 \bar{Q}(x,\Delta) &=& |a| + \delta \quad \mbox{if } \quad |x| > 1  \nonumber \\
 \bar{Q}(x,\Delta) &=& \alpha \quad \mbox{if } \quad 0 \leq |x| \leq 1, \Delta \geq L   \nonumber \\
 \bar{Q}(x,\Delta) &=& 1 \quad \quad  \mbox{if } \quad 0 \leq |x| \leq 1, \Delta < L, \nonumber
\end{eqnarray}

\begin{theorem}\label{InvNoiseless} \cite{YukMeynTAC2010} \cite{YukTAC2010}
Consider an adaptive quantizer applied to the linear control system described by (\ref{scalarSystem}). If the noiseless channel has capacity, \[R > \log_2(\lceil |a| \rceil + 1),\] and for the adaptive quantizer in (\ref{QuantizerUpdate2CHP7}), if the quantizer bin sizes are such that their (base$-2$) logarithms are integer multiples of some scalar $s$, and $\log_2(\bar{Q}(\cdot,\cdot))$ take values in integer multiples of $s$ where the integers taken are relatively prime (that is they share no common divisors except for $1$), then the process $\{(X_t,\Delta_t)\}$ is a positive (Harris) recurrent Markov chain (and has a unique invariant distribution).
\end{theorem}

In \cite{YukMeynTAC2010} it was shown that an $m$-small set (since a petite set in an irreducible and aperiodic Markov chain is $m$-small \cite{MeynBook}) can be constructed so that return conditions are satisfied. Hence, the return time properties directly leads to a stability result. The small set discussion in \cite{YukMeynTAC2010} builds on the Markovian property and irreducibility and aperiodicity of the Markov chain, together with a {\it uniform countable additivity} condition from \cite{Tweedie}.

We can obtain the stability result through the analysis in this paper, without defining a small/petite set: One can view the system as: $(\Delta_{t+1},x_{t+1})=F(\Delta_{t},x_{t},w_t)$, where the state is now $s_t:=(\Delta_{t},x_{t})$ and the independence of $w_t$ makes the process $(\Delta_{t},x_{t})$ Markov. Let the transition kernel be denoted with $P$. The finiteness of $\limsup_{t \to \infty} E[\Delta_t^2 + x_t^2]$ can be established by a Lyapunov analysis similar to \cite{YukTAC2010} and \cite{AndrewJohnstonReport}. However, $F$ here is not continuous in $s_t$.  Nonetheless, the set of discontinuity is given by: 
\[D = \bigg\{x, \Delta: {x \over \Delta} \in \{ -{K \over 2}, \cdots, {K \over 2} \}, \quad \Delta \in {\cal N} \bigg\}, \]
where ${\cal N}$ is the set of admissible bin sizes which is a countable set by the hypothesis of relative primeness. As a result $D$ is also countable and closed (since the elements are uniformly separated from each other). Furthermore, any weak limit of a converging sequence of expected occupational measures has zero measure on $D$, as can be deduced from the condition that every open set $D_{\epsilon} = \{x,\Delta: d((x, \Delta),D) < \epsilon\}$ is such that 
\begin{eqnarray}
v_{t_k}P (D_{\epsilon})  &&= \sum_{\Delta} \int_{z} v_{t_k}(dz,\Delta) \sum_{\Delta_{t_k}} P(\Delta_{t_k}|x_{t_k-1}=z, \Delta_{t_{k}-1}=\Delta) \nonumber \\
&& \quad \quad \quad \quad \times \sum_{m =-{K \over 2}}^{{K \over 2}} P( x_{t_k} \in [m\Delta_{t_k}-\epsilon,m\Delta_{t_k} + \epsilon]|x_{t_k-1}=z, \Delta_{t_{k}-1}=\Delta)  \nonumber \\
&& \leq L_1 \epsilon,  \nonumber 
\end{eqnarray}
for some $L_1 < \infty$ since $P(x_{t+1} \in dx| x, \Delta)$ has a density which is uniformly bounded for all $z, \Delta$ and the conditional probability $P(\Delta_{t_k}| x_{t_k-1}=z, \Delta_{t_{k}-1}=\Delta)$ has finite support.  By Theorem \ref{mainTheorem3}, the result follows. Finally, ergodicity follows from the irreducibility of the Markov process. \qed

\section{Conclusion}\label{Conclusion}

In this paper, a method to verify stochastic stability, asymptotic mean stationary and ergodicity properties of a class of non-Markovian stochastic processes has been introduced. Applications to practically important feedback coding schemes and networked control have been investigated. Further applications in control of non-Markovian systems and stability of non-linear filters are interesting research directions. 

\appendix

\section{Ergodic theorems for Markov chains}\label{appendixA}

Suppose that $\{X_t\}_{t\ge 0}$ denote a discrete-time Markov chain with state space
$\mathbb{X}$, a Polish space; its Borel $\sigma$-field is denoted by
${\cal B}(\mathbb{X})$, defined on a probability space $(\Omega, {\cal F}, {\bf P})$.  The transition probability is denoted by $P$, so that for any
$x$, $A \in {\cal B}(\mathbb{X})$, the probability of moving in one step from the state
$x$ to the set $A$ is given by $ {\bf P}(X_{t+1}\in A \mid x_t=x) = P(x,A)$.  The
$n$-step transitions are obtained via composition in the usual way, $
{\bf P}(X_{t+n}\in A \mid X_t=x) = P^n(x,A)$, for any $n\ge1$.  The transition law
acts on measurable functions $f\colon \mathbb{X}\to \mathbb{R}$ and measures $\mu$ on ${\cal B}(\mathbb{X})$
via  $Pf(x) := \int_{\mathbb{X}} P(x,dy) f(y),\quad x\in\mathbb{X},$ and $\mu P (A) := \int_{\mathbb{X}} \mu(dx)P(x,A)$, $A\in {\cal B}(\mathbb{X})$. A probability measure $\pi$ on ${\cal B}(\mathbb{X})$ is
called invariant if $\pi P= \pi$, i.e.,
\[
\int \pi(dx) P(x,A) = \pi(A),\qquad A\in {\cal B}(\mathbb{X}).
\]
For any initial probability measure $v$, by the Ionescu Tulcea theorem \cite{HernandezLermaMCP}, we can uniquely construct a stochastic
process with transition law $P$, and satisfying $X_0\sim v$.  We let $P_v$ denote the resulting probability measure on the sample space
$(\mathbb{X}, {\cal B}(\mathbb{X}))^{\mathbb{Z}_+}$, with the usual convention for $v=\delta_{x}$ when the initial state is $x\in \mathbb{X}$ in which case we write $P_x$ for the resulting probability measure. Likewise, $E_x$ denotes the expectation operator when the initial condition is given by $X_0=x$. When $v=\pi$,  the resulting process is stationary. 

When an invariant probability measure is known to exist for a Markov chain, we state the following ergodicity results.
 \begin{theorem}\label{convergenceT} \cite[Theorems 2.3.4-2.3.5]{HernandezLermaLasserre}
Let $\bar{P}$ be an invariant probability measure for a Markov process.
\begin{itemize}
\item[(i)] [Individual ergodic theorem] Let $X_0=x$. For every $f \in L_1(\bar{P})$ \[{1 \over N} E_x[\sum_{n=0}^{N-1} f(X_n)] \to f^*(x),\] 
for all $x \in B_f$ where $\bar{P}(B_f)=1$ (where $B_f$ denotes that the set of convergence may depend on $f$) for some $f^{*}$. 
\item[(ii)] [Mean ergodic theorem] Furthermore, the convergence ${1 \over N} E_x[\sum_{n=0}^{N-1} f(X_n)] \to f^*(x)$ is in $L_1(\bar{P})$.
\end{itemize}
\end{theorem}

 \begin{theorem}\label{convergenceT2} \cite[Theorem 2.5.1]{HernandezLermaLasserre}
Let $\bar{P}$ be an invariant probability measure for a Markov process. With $X_0=x$, fr every $f \in L_1(\bar{P})$ \[{1 \over N} \sum_{n=0}^{N-1} f(X_n) \to f^*(x),\] for all $x \in B_f$ where $\bar{P}(B_f)=1$ for some $f^{*}(x)$ with 
\[\int \bar{P}(dx) f^*(x) = \int \bar{P}(dx) f(x) \]
\end{theorem}

One may state further refinements; see \cite{HernandezLermaLasserre} for the locally compact case and \cite{worm2010ergodic} for the Polish state space case. 
 \begin{theorem}\label{convergenceT3} \cite{HernandezLermaLasserre} \cite{worm2010ergodic}
Let $\bar{P}$ be an invariant probability measure for a Markov process.
\begin{itemize}
\item[(i)] [Ergodic decomposition and weak convergence] For $x$, $\bar{P}$ a.s., ${1 \over N} E_x[\sum_{t=0}^{N-1} 1_{\{x_n \in \cdot\}}] \to P_x( \cdot)$ weakly and $\bar{P}$ is invariant for $P_x(\cdot)$ in the sense that
     \[\bar{P}(B) = \int P_x(B) \bar{P}(dx)   \]
\item[(ii)] [Convergence in total variation] For all $\mu \in {\cal P}(\mathbb{X}^{\mathbb{N}})$ which satisfies that $\mu \ll \bar{P}$ (that is, $\mu$ is
 absolutely continuous with respect to the stationary mean), there exists $v^*$ such that
\[ \| E_{\mu}[{1 \over N}\sum_{t=0}^{N-1} 1_{\{T^nX \in \cdot\}}]  - v^*(\cdot) \|_{TV} \to 0. \]
\end{itemize}
\end{theorem}


%


\end{document}